\documentclass[12pt]{iopart}

\usepackage{iopams}
\usepackage{cite}
\usepackage{amsthm}
\usepackage{enumerate}
\newtheorem{theorem}{Theorem}

\newtheorem{exmp}{Example}
\newtheorem{defn}{Definition}

\theoremstyle{definition}

\begin{document}

\title[Exact solutions of laplace equations by DTM]{Exact solutions of Laplace equation by differential transform method}

\author{M. Jamil Amir, M. Yaseen and Rabia Iqbal}

\address{Department of Mathematics, University of Sargodha,  Sargodha 40100, Pakistan.}
\ead{mjamil.dgk@gmail.com}
\begin{abstract}
In this paper, we solve Laplace equation analytically by using differential transform method. For this purpose, we consider four models with two Dirichlet and two Neumann boundary conditions and obtain the corresponding exact solutions. The obtained results show the simplicity of the method and massive reduction in calculations when one compares it with other iterative methods, available in literature. It is worth mentioning that here only a few number of iterations are required to reach the closed form solutions as series expansions of some known functions.
\end{abstract}

\maketitle
\section{Introduction}
Most of the problems in physics and engineering are usually expressed in the form of partial differential equations which are often too complicated to be solved analytically. Even if one may obtain the exact solution of some problems it involves much tedious calculations or it is difficult to interpret the solutions. The struggle to overcome this difficulty led to the invention of different approximate analytical methods, such as, Adomian decomposition  method (ADM) \cite{adomian}, homotopy analysis method (HAM) \cite{Liu1}, variational iteration method (VIM) \cite{vim1},
the new iterative method (NIM) \cite{varsha1, varsha2}  and differential transform method (DTM) \cite{zhou} etc.

The DTM was first introduced by Zhou\cite{zhou} in 1986 for solving linear and non-linear initial value problems in electrical circuit analysis. The method has been extensively used by researchers to solve linear and non-linear ordinary differential equations \cite{siddiqi,chen,chiou,chen1,Kuo,Kuo1}. In 1999, Chen an Ho \cite{chen2} extended this method to two dimensional differential transform method for solving the partial differential and integral equations. This method has been used for solving linear and non linear Goursat problem \cite{Taghvafard}, for linear partial differential equations of factorial order \cite{Odibat}, for solving system of differential equations \cite{Ayaz} and for differential-algebraic equations \cite{Ayaz1}.

The Laplace equation is the best model to study most of the problems of different branches of physics and mechanics, for example, heat and mass transfer, fluid mechanics, elasticity and electrostatics. The two dimensional Laplace equation
$$\nabla^2 u=0,$$ where $\nabla^2=\frac{\partial^2 u}{\partial x^2}+\frac{\partial^2 u}{\partial y^2}$ is Laplacian, is usually considered with Dirichlet and Neumann boundary conditions. The Laplace equation has already been examined using several iterative methods such as the new iterative method \cite{Yaseen}, homotopy analysis method \cite{INC}, Adomian decomposition method \cite{sadighi} and variational iteration method \cite{wazwaz}.

The purpose of present paper is to apply two dimensional DTM to four models of  Laplace equation, two with Dirichlet and two with Neumann boundary conditions. In each case, this method constructs one analytic solution  without requiring liberalization or discretization. Another important advantage of DTM is simplicity in its algorithm. Moreover, it greatly reduces the size of computational work compared with the existing iterative methods and provide accurate series solution with rapid convergence. In next section, we present two dimensional differential transform method \cite{Ming}.
 \section{Two dimensional differential transform method}
Consider a function of two variables $u(x,y)$ and suppose that it can be represented as a product of two
single-variable functions, i.e., $u(x,y)=f(x)g(y)$. Based on the properties of One-dimensional differential transform,
function $u(x,y)$ can be represented as
\begin{equation}
u(x,y)=\sum\limits_{m=0}^{\infty}F(m)x^{m}\sum\limits_{n=0}^{\infty}G(n)y^{n}=\sum\limits_{m=0}^{\infty}\sum\limits_{n=0}^{\infty}U(m,n)x^{m}y^{n}
\end{equation}
where $U(m,n)=F(m)G(n)$ is called spectrum of $u(x,y)$. \\
 The basic definitions and operations of two-dimensional differential transform are introduced as follows:
\begin{defn}
 If a function $u(x,y)$ is analytic and differentiated continuously with respect to time $t$ in the domain of interest and let
\begin{equation}\label{1}
   U(m,n)=\frac{1}{m!n!}\left[\frac{\partial^{m+n}u(x,y)}{\partial x^{m}\partial y^{n}}_{x=x_{0},y=y_{0}}\right],
\end{equation}
where the spectrum $U(m,n)$ is the transformed function. Then the differential inverse transform of $U(m,n)$ is defined as follows
\begin{equation}\label{2.2}
u(x,y)=\sum\limits_{m=0}^{\infty}\sum\limits_{n=0}^{\infty}U(m,n)(x-x_{0})^{m}(y-y_{0})^{n}
\end{equation}
\end{defn}
The lower case $u(x,y)$ represents the original function while the uppercase $U(m,n)$ stands for the transformed function. Combining Eq.(\ref{1}) and Eq.(\ref{2.2}), it can be obtained that
\begin{equation}  \label{2.3}
u(x,y)=\sum\limits_{m=0}^{\infty}\sum\limits_{n=0}^{\infty}\frac{1}{m!n!}\left[\frac{\partial^{m+n}u(x,y)}{\partial x^{m}\partial y^{n}}_{x=x_{0},y=y_{0}}\right](x-x_{0})^{m}(y-y_{0})^{n}
\end{equation}
On the basis of above definition we have following fundamental results of two dimensional differential transform (see \cite{Ming})
\begin{theorem}  \label{2.1}
Suppose that $U(m,n)$, $V(m,n)$ and $W(m,n)$ are the differential transforms of the functions $u(x,y)$, $v(x,y)$ and $w(x,y)$ respectively at $(0,0)$, then it follows that
\begin{enumerate}[(a)]
\item If $u(x,y)=v(x,y)\pm w(x,y)$, then $$U(m,n)=V(m,n)\pm W(m,n)$$\\
\item If $u(x,y)=av(x,y)$, then $$U(m,n)=aV(m,n)$$\\
\item If $u(x,y)=v(x,y)w(x,y)$, then
$$U(m,n)=\sum\limits_{k=0}^{m}\sum\limits_{l=0}^{n}V(k,n-l)W(m-k,l)$$\\
\item If $u(x,y)=\frac{\partial^{r+s}v(x,y)}{\partial x^{r}\partial y^{s}}$, then \\
$$U(m,n)=\frac{(m+r)!}{m!}\frac{(n+s)!}{n!}V(m+r,n+s)$$\\
\item If $u(x,y)=e^{av(x,y)}$, then\\
\begin{eqnarray*}
U(m,n)=\left\{\begin{array}{ll}
                e^{av(0,0)}, & m=n=0 \\
                a\sum\limits_{k=0}^{m-1}\sum\limits_{l=0}^{n}\frac{m-k}{m}V(m-k,l)U(k,n-l), & m\geq1 \\
                a\sum\limits_{k=0}^{m}\sum\limits_{l=0}^{n-1}\frac{n-l}{n}V(k,n-l)U(m-k,n), & n \geq 1.
\\
              \end{array}\right.
\end{eqnarray*}
\item If $u(x,y)=x^{k}y^{h}$, then\\
\begin{eqnarray*}
U(m,n)=\left\{\begin{array}{cc}
                \delta(m-k,n-h), &  m=k, n=h, \\
                 0,              & otherwise \\

               \end{array}\right.
\end{eqnarray*}
\item If $u(x,y)=x^{k}e^{ay}$, then\\ $$U(m,n)=\delta(m-k)\frac{a^{n}}{n!}$$
\end{enumerate}
\end{theorem}
\section{Solutions of Laplace Equation}
In this section we apply DTM to four physical models of Laplace equation to establish exact solutions.
\begin{exmp}
Consider the second order Laplace equation, given as
\begin{equation}
\ \ u_{xx}+u_{yy}=0,~~~~~~0<x,y<\pi  \label{3.1}
\end{equation}
with Dirichlet  boundary conditions
\begin{eqnarray}\label{3.2}
\begin{array}{l}
u(x,0)=\sinh x,~u(x,\pi)=-\sinh x,\\
u(0,y)=0,~~~u(\pi,y)=\sinh(\pi)\cos y.
\end{array}
\end{eqnarray}
\end{exmp}
Taking the differential transform of (\ref{3.1}) and using Theorem 2.1, it can be obtained that
\begin{equation}\label{3.3}
(m+1)(m+2)U(m+2,n)+(n+1)(n+2)U(m,n+2)=0.
\end{equation}
The conditions (\ref{3.2}) and Eq. (\ref{2.2}) imply that
\begin{equation}
u(x,0)=\sum\limits_{m=0}^{\infty}U(m,0)x^{m}=\sinh x=\sum\limits_{m=1,3,5\cdots}^{\infty}\frac{x^m}{m!},
\end{equation}
which, on comparing the both sides yields
\begin{eqnarray}\label{3.4}
U(m,0)=\left\{\begin{array}{cc}
                    \frac{1}{m!} & when~m~is~odd \\
                    0 & otherwise
                  \end{array}\right.
\end{eqnarray}
Also Eq. (\ref{3.2}) and Eq. (\ref{2.2}) imply that
\begin{equation}
u(0,y)=\sum\limits_{n=0}^{\infty}U(0,n)y^{n}=0,
\end{equation}
which yields
\begin{equation}\label{3.5}
U(0,n)=0.
\end{equation}
Substituting Eqs.(\ref{3.4}) and (\ref{3.5}) into Eq. (\ref{3.3}) and after some calculations, we reach at
\begin{eqnarray}\label{3.6}
U(m,n)=\left\{\begin{array}{cc}
                    \frac{(-1)^{\frac{n}{2}}}{m!n!}, & when~m~is~odd~and~n~is~even \\
                    0, & otherwise.
                  \end{array}\right.
\end{eqnarray}
Substituting Eq.(\ref{3.6}) into Eq.(\ref{2.2}), we have
\begin{eqnarray}  \label{3.7}
u(x,y)&=&\sum\limits_{m =1,3,5,\cdots}^{\infty}\sum\limits_{n=0,2,4,\cdots}^{\infty}\frac{(-1)^{\frac{n}{2}}}{m!n!}x^{m}y^{n}, \nonumber \\
&=&\left(\sum\limits_{m=1,3,5\cdots}^{\infty}\frac{x^{m}}{m!}\right)\left(\sum\limits_{n=0,2,4\cdots}^{\infty}\frac{(-1)^{\frac{n}{2}}}{n!}y^{n}\right),\nonumber \\
&=&\sinh x\cos y,
\end{eqnarray}
which is the required solution.
\begin{exmp}
Consider the second order Laplace equation,
\begin{equation}
\ \ u_{xx}+u_{yy}=0,~~~~~~0<x,y<\pi  \label{4.1}
\end{equation}%
with Dirichlet boundary conditions
\begin{eqnarray}\label{4.2}
\begin{array}{l}
u(x,0)=0,~~~u(x,\pi)=0,\\
u(0,y)=\sin y,u(\pi,y)=\cosh(\pi)\sin y.
\end{array}
\end{eqnarray}
\end{exmp}
Taking the differential transform of (\ref{4.1}) and using Theorem 2.1, it follows that
\begin{equation}\label{4.3}
(m+1)(m+2)U(m+2,n)+(n+1)(n+2)U(m,n+2)=0.
\end{equation}
From conditions (\ref{4.2}) and Eq. (\ref{2.2}), we obtain
\begin{equation}
u(x,0)=\sum\limits_{m=0}^{\infty}U(m,0)x^{m}=0,
\end{equation}
which, on comparing the both sides yields 

\begin{equation}\label{4.4}
U(m,0)=0.
\end{equation}
Also, Eqs. (\ref{2.2}) and (\ref{4.2}), imply that
\begin{equation}
u(0,y)=\sum\limits_{n=0}^{\infty}U(0,n)y^{n}=\sin y=\sum \limits_{0}^{\infty}\frac{(-1)^{\frac{n-1}{2}}}{n!}y^{n},
\end{equation}
which on comparing the both sides results as
\begin{eqnarray}\label{4.5}
U(0,n)=\left\{\begin{array}{cc}
                    \frac{(-1)^{\frac{n-1}{2}}}{n!}, & when~n~is~odd \\
                    0, & otherwise.
                  \end{array}\right.
\end{eqnarray}
Substituting Eqs. (\ref{4.4}) and (\ref{4.5}) in Eq. (\ref{4.3}) and making some calculations, we have
\begin{eqnarray}\label{4.6}
U(m,n)=\left\{\begin{array}{cc}
                    \frac{(-1)^{\frac{n-1}{2}}}{m!n!}, & when~m~is~even~and~n~is~odd \\
                    0, & otherwise.
                  \end{array}\right.
\end{eqnarray}
Now making use of Eq. (\ref{4.6}) in Eq. (\ref{2.2}), we obtain
\begin{eqnarray}  \label{3.7}
u(x,y)&=&\sum\limits_{m=0,2,4\cdots}^{\infty}\sum\limits_{n=1,3,5\cdots}^{\infty}\frac{(-1)^{\frac{n-1}{2}}}{m!n!}x^{m}y^{n}, \nonumber \\
&=&\left(\sum\limits_{m=0,2,4\cdots}^{\infty}\frac{x^{m}}{m!}\right)\left(\sum\limits_{n=1,3,5\cdots}^{\infty}\frac{(-1)^{\frac{n-1}{2}}}{n!}y^{n}\right), \nonumber \\
&=&\cosh x\sin y,
\end{eqnarray}
which is the exact solution.
\begin{exmp}
Consider the second order Laplace equation,
\begin{equation}
\ \ u_{xx}+u_{yy}=0,~~~~0<x,y<\pi  \label{5.1}
\end{equation}%
with Neumann boundary conditions
\begin{eqnarray}\label{5.2}
\begin{array}{l}
u_y(x,0)=0,~u_y(x,\pi)=2\cos 2x\sinh 2\pi\\
u_x(0,y)=0,~u_x(\pi,y)=0
\end{array}
\end{eqnarray}
\end{exmp}

Taking the differential transform of (\ref{5.1}) and using Theorem 2.1, it follows that
\begin{equation}\label{5.3}
(m+1)(m+2)U(m+2,n)+(n+1)(n+2)U(m,n+2)=0.
\end{equation}
Eqs. (\ref{2.2}) and  (\ref{5.2}) imply that
\begin{eqnarray}
u_{y}(x,\pi)&=&\sum\limits_{m=0}^{\infty}\sum\limits_{n=0}^{\infty}n\pi^{n-1}U(m,n)x^{m}, \nonumber \\
&=&2\cos2x\sinh2\pi, \nonumber \\ &=&\sum\limits_{m=0}^{\infty}\frac{(-1)^{\frac{m}{2}}(2x)^m}{m!}\sum\limits_{n=0}^{\infty}\frac{(2\pi)^n}{n!}.
\end{eqnarray}
By changing the index $n$ , and comparison,  we have
$$U(m,n+1)=\frac{(-1)^{\frac{m}{2}}2^{{m+n+1}}}{(n+1)m!n!},$$
which yields to following general relation
\begin{eqnarray}\label{5.4}
U(m,n)=\left\{\begin{array}{cc}
                    \frac{(-1)^{\frac{m}{2}}2^{{m+n}}}{m!n!}, & when~m~and~n~are~even \\
                    0, & otherwise.
                  \end{array}\right.
\end{eqnarray}

Substituting Eq. (\ref{5.4}) in Eq. (\ref{2.2}), we obtain
\begin{eqnarray}\label{5.5}
u(x,y)&=&\sum\limits_{m=0,2,4,\cdots}^{\infty}\sum\limits_{n=0,2,4,\cdots}^{\infty}\frac{(-1)^{\frac{m}{2}}2^{m+n}}{m!n!}x^{m}y^{n}, \nonumber \\
&=&\left(\sum\limits_{m=0,2,4\cdots}^{\infty}\frac{(-1)^{\frac{m}{2}}(2x)^{m}}{m!}\right)\left(\sum\limits_{n=0,2,4\cdots}^{\infty}\frac{(2y)^{n}}{n!}\right), \nonumber \\
&=&\cos2x\cosh 2y,
\end{eqnarray}
which is the exact solution.
\begin{exmp}
Consider the second order Laplace equation,
\begin{equation}
\ \ u_{xx}+u_{yy}=0,~~~~~0<x,y<\pi  \label{6.1}
\end{equation}%
with Neumann boundary conditions
\begin{eqnarray}\label{6.2}
\begin{array}{l}
u_y(x,0)=\cos x,~u_y(x,\pi)=\cosh \pi\cos x,\\
u_x(0,y)=0,~~~~u_x(\pi,y)=0.
\end{array}
\end{eqnarray}
\end{exmp}
Taking the differential transform of (\ref{6.1}) and using Theorem 2.1, it yields
\begin{equation}\label{6.3}
(m+1)(m+2)U(m+2,n)+(n+1)(n+2)U(m,n+2)=0.
\end{equation}
From Eqs. (\ref{2.2}) and (\ref{6.2}), we obtain that
\begin{equation}
u_{y}(x,0)=\sum\limits_{m=0}^{\infty}U(m,1)x^{m}=\cos x=\sum\limits_{m=0}^{\infty}\frac{(-1)^{\frac{m}{2}}x^m}{m!},
\end{equation}
which on comparison with the cosine series, yields
\begin{eqnarray}\label{6.4}
U(m,1)=\left\{\begin{array}{cc}
                    \frac{(-1)^{\frac{m}{2}}}{m!}, & when~m~is~even \\
                    0, & otherwise.
                  \end{array}\right.
\end{eqnarray}
Also, Eqs. (\ref{2.2}) and (\ref{6.2}) imply that
\begin{equation}
u_{x}(0,y)=\sum\limits_{n=0}^{\infty}U(1,n)y^{n}=0,
\end{equation}
which yields
\begin{equation}\label{6.5}
U(1,n)=0.
\end{equation}
Substituting Eqs. (\ref{6.4}) and (\ref{6.5}) in Eq. (\ref{6.3}), we obtain the general relation as
\begin{eqnarray}\label{6.6}
U(m,n)=\left\{\begin{array}{cc}
                    \frac{(-1)^{\frac{m}{2}}}{m!n!}, & when~m~is~even~and~n~is~odd \\
                    0, & otherwise.
                  \end{array}\right.
\end{eqnarray}
making use of Eq. (\ref{6.6}) in Eq. (\ref{2.2}), we get
\begin{eqnarray}  \label{6.7}
u(x,y)&=&\sum\limits_{m=0,2,4\cdots}^{\infty}\sum\limits_{n=1,3,5,...}^{\infty}\frac{(-1)^{\frac{m}{2}}}{m!n!}x^{m}y^{n}, \nonumber \\
&=&\left(\sum\limits_{m=0,2,4\cdots}^{\infty}\frac{(-1)^{\frac{m}{2}}x^{m}}{m!}\right)\left(\sum\limits_{n=1,3,5,\cdots}^{\infty}\frac{y^{n}}{n!}\right), \nonumber \\
&=&\cos x\sinh y,
\end{eqnarray}
which is the exact solution.
\section{Conclusion}
In this paper, we have solved the four cases of Laplace equation which have been often used in describing the physical phenomenon of real life. We have successfully developed the DTM to obtain the exact solutions of Laplace equation. The method gives rapid convergence by using minimum number of iterations. One extra advantage is massive reduction of calculation compared with other iterative methods, available in literature. It is apparent that this method is very powerful and efficient for solving different kinds of physical problems that arise in physics and engineering.

\end{document}